\documentclass[12pt]{extarticle}

\usepackage[margin=1.5in]{geometry}  
\usepackage{graphicx}              
\usepackage{amsmath}               
\usepackage{amsfonts}              
\usepackage{amsthm}                
\usepackage{framed}

\begin{document}

\nocite{*}

\title{\bf An Elementary Proof of an Asymptotic Formula of Ramanujan}

\author{\textsc{Adrian W. Dudek} \\ 
Mathematical Sciences Institute \\
The Australian National University \\ 
\texttt{adrian.dudek@anu.edu.au}}
\date{}

\maketitle

\begin{abstract}
We derive the mean square of the divisor function using only elementary techniques.
\end{abstract}

\section{Introduction}

Let $d(n)$ denote the number of divisors of a positive integer $n$. In 1915, Ramanujan \cite{ramanujan} stated (without proof) the formula

\begin{equation} \label{main}
\sum_{n \leq x} d^2(n) = A x \log^3 x + Bx \log^2 x + C x \log x+ Dx +O(x^{3/5 + \epsilon})
\end{equation}

where $A = \pi^{-2}$, and $B$, $C$ and $D$ are more complicated constants (expressible, nonetheless). In particular, the \textit{average value} of the square of the divisor function is $\pi^{-2}  \log^3 n.$ It is the point of this note to prove this average value using only elementary techniques.

The formula (\ref{main}) has been studied by others (Wilson \cite{wilson}, Suryanarayana \& Sita Rama Chandra Rao \cite{suryanarayana}) for the problem of reducing the order of the error term, akin somewhat to Dirichlet's divisor problem. In each case, the formula was derived again independently. Wilson used analytic techniques, namely Perron's formula and Cauchy's integral theorem for his result. Suryanarayana and Sita Rama Chandra Rao started with the convolution formula

\begin{equation} \label{convolution}
d^2 (n) = \sum_{\delta^2 | n} \mu(\delta) d_4 \Big(\frac{n}{\delta^2}\Big),
\end{equation}

where $\mu(n)$ is the Moebius function and $d_k(n)$ denotes the number of ways one can write $n$ as the product of $k$ positive integers, and used, among other things, analytic bounds for the Mertens function given by Landau \cite{landau}. 

\section{An Elementary Proof}

We start with $(\ref{convolution})$ (see Lemma 3.1 of \cite{suryanarayana} for an elementary proof) and use only the two elementary bounds

\begin{eqnarray} \label{harmonic}
\sum_{n \leq x} \frac{\log^k n}{n} & = & \frac{\log^{k+1} x}{k+1} +O(1)
\end{eqnarray}

which holds for all $k \geq 0$,

\begin{eqnarray} \label{divisor}
\sum_{n \leq x} d(n) & = & x \log x + O(x), \label{divisor} 
\end{eqnarray}

and the fact that $\zeta(2) = \pi^2/6$. A proof of $(\ref{harmonic})$ follows from comparing the sum to the corresponding integral, whereas a proof of $(\ref{divisor})$ can be found in most number theory textbooks and usually forms a preface to Dirichlet's hyperbola method. In particular, one can see Lemma 3.3 in Apostol's book \cite{apostol}.  There also exists a good stash of elementary proofs which give the exact value of $\zeta(2)$ (see the recent proof by Daners \cite{daners}, for example).

We sum each side of (\ref{convolution}) to construct the relation

\begin{eqnarray} \label{d}
\sum_{n \leq x} d^2(n) & = & \sum_{m \delta^2 \leq x} \mu(\delta) d_4 (m) \nonumber \\
& = & \sum_{\delta \leq \sqrt{x}} \mu(\delta) \sum_{m \leq x/\delta^2} d_4 (m).
\end{eqnarray}

Our plan of attack was born of the realisation that the inner sum can be connected to a sum of $d_2 (n)=d(n)$ by the classic observation that

$$\sum_{n \leq y} d_k (n) = \sum_{mn \leq y} d_{k-1}(n),$$

where the sum on the right hand side is over both $m$ and $n$. We can apply this observation twice to the inner sum in $(\ref{d})$ to get

\begin{eqnarray} \label{bum}
\sum_{m \leq x/\delta^2} d_4 (m) & = & \sum_{qm \leq x/\delta^2} d_3(m) \nonumber \\ 
& = & \sum_{q \leq x/\delta^2} \sum_{m \leq x/(\delta^2 q)} d_3(m) \nonumber \\
& = &  \sum_{q \leq x/\delta^2} \sum_{r \leq x/(\delta^2 q)} \sum_{m \leq x/(\delta^2 q r)} d(m).
\end{eqnarray}

We estimate $(\ref{bum})$, by working from the innermost sum outwards, before placing our result into $(\ref{d})$. By $(\ref{divisor})$ we have that

$$\sum_{m \leq x/(\delta^2 q r)} d(m) = \frac{x}{\delta^2 q r} \log\Big( \frac{x}{\delta^2 q r}\Big) + O\Big( \frac{x}{\delta^2 q r} \Big).$$

As the innermost sum is now estimated, we move to the next. We have

\begin{eqnarray*}
 \sum_{r \leq x/(\delta^2 q)} \sum_{m \leq x/(\delta^2 q r)} d(m) & = &  \sum_{r \leq x/(\delta^2 q)} \bigg[ \frac{x}{\delta^2 q r} \log\Big( \frac{x}{\delta^2 q r}\Big) + O\Big( \frac{x}{\delta^2 q r} \Big)\bigg].
\end{eqnarray*}

By writing 

$$\log\Big( \frac{x}{\delta^2 q r} \Big) = \log\Big( \frac{x}{\delta^2 q}\Big) - \log r$$

we thus have

\begin{eqnarray*}
 \sum_{r \leq x/(\delta^2 q)} \sum_{m \leq x/(\delta^2 q r)} d(m) & = & \frac{x}{\delta^2 q} \log\Big( \frac{x}{\delta^2 q} \Big) \sum_{r \leq x/(\delta^2 q)} \frac{1}{r} - \frac{x}{\delta^2 q} \sum_{r \leq x/(\delta^2 q)} \frac{\log r}{r} \\ \\
 & + & O\Big( \sum_{ r \leq x/(\delta^2 q)} \frac{x}{\delta^2 q r} \Big).
\end{eqnarray*}

Each of these sums can be estimates using $(\ref{harmonic})$ with either $k=0$ (the first and third sum) or $k=1$ (the second sum) to arrive at

\begin{eqnarray*}
 \sum_{r \leq x/(\delta^2 q)} \sum_{m \leq x/(\delta^2 q r)} d(m) & = & \frac{1}{2} \frac{x}{\delta^2 q} \log^2 \Big( \frac{x}{\delta^2 q} \Big) + O\Big( \frac{x}{\delta^2 q} \log\Big( \frac{x}{\delta^2 q} \Big) \Big).
\end{eqnarray*}

We will complete our evaluation of $(\ref{bum})$ now with

\begin{eqnarray} \label{bumbo}
\sum_{m \leq x/\delta^2} d_4 (m) & = & \sum_{q \leq x/\delta^2} \sum_{r \leq x/(\delta^2 q)} \sum_{m \leq x/(\delta^2 q r)} d(m) \nonumber \\
& = & \sum_{q \leq x/\delta^2} \bigg[  \frac{1}{2} \frac{x}{\delta^2 q} \log^2 \Big( \frac{x}{\delta^2 q} \Big) + O\Big( \frac{x}{\delta^2 q} \log\Big( \frac{x}{\delta^2 q} \Big) \Big) \bigg].
\end{eqnarray}

This works in much the same way as before, however this time we write

$$\log^2\Big( \frac{x}{\delta^2 q} \Big) = \log^2 \Big( \frac{x}{\delta^2} \Big) + \log^2 q - 2 \log \Big( \frac{x}{\delta^2} \Big) \log q.$$

Inserting this into (\ref{bumbo}) and estimating using (\ref{harmonic}) with $k=0, 1$ and $2$ yields

\begin{eqnarray*} 
\sum_{m \leq x/\delta^2} d_4 (m) & = & \frac{1}{6} \frac{x}{\delta^2} \log^3 \Big( \frac{x}{\delta^2} \Big) + O\Big( \frac{x}{\delta^2} \log^2 \Big( \frac{x}{\delta^2} \Big)\Big).
\end{eqnarray*}

Finally, by (\ref{d}) we have

\begin{eqnarray*} 
\sum_{n \leq x} d^2(n) & = &  \sum_{\delta \leq \sqrt{x}} \mu(\delta) \bigg[ \frac{1}{6} \frac{x}{\delta^2} \log^3 \Big( \frac{x}{\delta^2} \Big) + O\Big( \frac{x}{\delta^2} \log^2 \Big( \frac{x}{\delta^2} \Big)\Big) \bigg] \\
& = & \frac{1}{6} x \log^3 x \sum_{\delta \leq \sqrt{x}} \frac{\mu(n)}{n^2} + O(x \log^2 x).
\end{eqnarray*}

It's also known that

\begin{eqnarray*}
\sum_{\delta \leq \sqrt{x}} \frac{\mu(n)}{n^2} & = & \frac{1}{\zeta(2)} + O\Big( \int_{\sqrt{x}}^{\infty} \frac{dt}{t^2} \Big) \\
& = & \frac{6}{\pi^2} + O\Big( \frac{1}{\sqrt{x}} \Big). 
\end{eqnarray*}

Substituting this in gives us our result:

$$\sum_{n \leq x} d^2(n) = \frac{1}{\pi^2} x \log^3 x + O(x \log^2 x).$$

As an aside, it should be possible to replace (\ref{divisor}) with the more accurate estimate

$$\sum_{n \leq x} d(n) = x \log x + (2 \gamma - 1) +O(x^{1/2})$$

to furnish the full asymptotic expansion. We leave such a venture to the more enthusiastic readers.

\bibliographystyle{plain}

\bibliography{biblio}

\begin{thebibliography}{1}

\bibitem{apostol}
T.~M. Apostol.
\newblock {\em Introduction to {A}nalytic {N}umber {T}heory}.
\newblock Springer, 1976.

\bibitem{daners}
D.~Daners.
\newblock A short elementary proof of $\sum 1/k^2 = \pi^2/6$.
\newblock {\em Mathematics Magazine}, 85(5):361--364, 2012.

\bibitem{landau}
E.~Landau.
\newblock Handbuch der {L}ehre von der {V}erteilung der {P}rimzahlen.
\newblock {\em Leipzig}, 1, 1909 (Chelsea reprint, 1953).

\bibitem{ramanujan}
S.~Ramanujan.
\newblock Some formulae in the analytic theory of numbers.
\newblock {\em Messenger of Math.}, 45:81--84, 1915.

\bibitem{suryanarayana}
D.~Suryanarayana and R.~Sita Rama Chandra~Rao.
\newblock On an asymptotic formula of {R}amanujan.
\newblock {\em Math. Scand}, 82:258--264, 1978.

\bibitem{wilson}
B.~M. Wilson.
\newblock Proofs of some formulae enunciated by {R}amanujan.
\newblock {\em Proc. Lond. Math. Soc.}, 2(21):235--255, 1922.

\end{thebibliography}

\end{document}